\def\Z{\mathbb Z}

    \def\N{\mathbb N}

  \def\G{\Gamma}

    \documentclass[a4paper]{article}
            \usepackage{amsfonts,amssymb}
\usepackage{amsmath}

\newtheorem{thm}{Theorem}
  
  \newtheorem{lem}[thm]{Lemma}

\begin{document}

  \title{On the finite images of finitely generated perfect groups}
  \author{Nikolay Nikolov }
    \date{}
  \maketitle

 \begin{abstract} Let $ d \geq 2$ be an integer. We conjecture that there is a finitely generated perfect group whose images include all finite $d$-generated perfect groups. We prove a special case of this conjecture for the finite perfect groups with a subnormal series of bounded length and factors which are abelian or semisimple.

 \end{abstract}
 Let $d \geq 2$ be an integer.
We say that group $G$ is perfect if $G$ equals its derived subgroup $[G,G]$. 

\textbf{Conjecture A} \emph{Let $d \geq 2$ be an integer. There is a finitely generated perfect group $\G _d$ whose homomorphic images include all $d$-generated finite perfect groups.}

One can interpret this conjecture in the theory of profinite groups as follows. A finitely generated profinite group $G$ is perfect if and only if there exists $d \in \mathbb N$ such that $G$ is an inverse limit of $d$-generated finite perfect groups.

\textbf{Conjecture B1} \emph{For every integer $d>1$ there is a finitely generated perfect profinite group whose finite images include all finite perfect $d$-generated groups.}

\textbf{Conjecture B2} \emph{Every finitely generated perfect profinite group contains a dense finitely generated perfect subgroup.}

We note that Conjecture A is equivalent to the simultaneous validity of Conjectures B1 and B2.

It is immediate that Conjectures B1 and B2 together imply Conjecture A. Conversely Conjecture A implies Conjecture B1 trivially: the required profinite group is the profinite completion $ \hat \G_d$ of $\G_d$. Finally Conjecture A implies Conjecture B2 as well, since a standard compactness argument gives that any $d$-generated perfect profinite group is a homomorphic image of $\hat \G_d$ which in turn contains an image of $\G_d$ as a dense perfect subgroup.

We remark that it is an open question of Andrei Jaikin-Zapirain \cite{problems} whether a non-abelian free profinite group is locally indicable. This will be answered negatively if Conjecture 2B is true for virtualy pro-nilpotent perfect profinite groups. Indeed, the universal Frattini cover $\mathcal F_G$ of a perfect finite group $G$ is a finitely generated, perfect, projective profinite group, see Theorem 1.1 in \cite{ef}. Therefore $\mathcal F_G$ is a closed subgroup of a free profinite group and Conjecture 2B provides a finitely generated perfect subgroup of $\mathcal F_G$.

While Conjecture 2B remains open, Theorem \ref{main} below provides an affirmative answer in the special case of perfect poly-(abelian or semisimple) profinite groups.

The first test of Conjecture A is to verify it for the family of non-abelian finite simple groups. This is a result of D. Segal.
\begin{thm}[\cite{segal}] \label{seg}
    There is a perfect group generated by 61 elements which has all finite non-abelian simple groups among its quotients.
\end{thm}

In this note we generalise Theorem \cite{segal} to the family of finite $d$-generated perfect groups which are obtained by a bounded number of successive extensions of abelian or semisimple groups. The precise definition is as follows.

For a finite group $G$ let \[ G_*:= \cap \{ N \triangleleft G \ |  \ G/N \textrm{ is abelian or simple} \}.\] Thus $G/G_*= A \times S$, where $A$ is isomorphic to the abelianisation $G/[G,G]$ of $G$ and $S$ is the largest semisimple quotient of $G$. Define the series of characteristic subgroups $(G_n)$ of $G$ by $G_0=G$ and $G_{n+1}=(G_n)_*$ for all integers $n \geq 0$. 

Let $\mathcal Y(d,k)$ be the set of (isomorphism classes) of finite $d$-generated perfect groups $G$ such that $G_k=\{1\}$. 

\begin{thm} \label{main}
    Given integers $k \geq 0,d \geq 1$ there is a finitely generated perfect group $\G_{d,k}$ whose finite images include $\mathcal Y(d,k)$.
\end{thm}

We will need several results which we list below.

\begin{lem}[Gaschutz \cite{Gas}]\label{gaschutz} Let $G$ be a finite group with a normal subgroup $N$. Assume that $G$ has a generating set of size $d \in \N$. Let $k \geq d$ be an integer and let $a_1N, \ldots, a_kN$ be a set of generators for $G/N$. There exist $b_i \in a_iN$ such that $b_1, \ldots b_k$ is a set of generators for $G$.    \end{lem}

\begin{lem} \label{G-1} Let $G$ be a group with generating set $\{ g_1, \ldots, g_d\}$. Let $M$ be a left $\Z G$-module. 

(a) $M(G-1)= \sum_{i=1}^d M(g_i-1)$

(b) If $M$ is generated by $m_1, \ldots, m_k$ as a $\Z G$-module then $M(G-1)$ is generated as $\Z G$-module by the elements $m_j(g_i-1)$ for $j=1,\ldots, k, i=1, \ldots, d$.

(c) If $G$ is a perfect group then $V=M(G-1)$ is a perfect $\Z G$-module, i.e. $V=V(G-1)$.
\end{lem}

\begin{thm} [\cite{LS}] \label{cover} There is constant $C$ such that for every non-abelian finite simple group $S$ and every normal subset $X \subseteq S$ the covering number of $S$ with respect to $X$ is at most $C \log |S|/ \log |X|$.
    
\end{thm}
The proof of Theorem \ref{main} is by induction on $k$. The case $k=0$ is clear since $\mathcal Y(d,0)$ is just the trivial group.
Assume that $k \geq 1$ and the Theorem holds for $k-1$. Let $\{G_j \ | \  j \in \N\} = \mathcal Y(d,k)$. For each $j \in \N$ let $W_j=A_j \times S_j$ be a normal subgroup of $G_j$ such that $G_j/W_j \in \mathcal Y(d,k-1)$ with $A_j$ abelian and $S_j$ semisimple. Note that if $k=1$ then $A_j=\{1\}$ since in that case $W_j=G_j$ is a perfect group. 

By the induction hypothesis there is a finitely generated perfect group $\G_{d,k-1}$ which maps onto each $G_j/W_j$.

We claim that without loss of generality we may assume that the group $\G_{d,k-1}$ is finitely presented with equal number of generators and relations.

Indeed, suppose $\G_{d,k-1}$ has generators $g_1, \ldots, g_m$. Let $F$ be the free abstract group on $m$ generators $x_1, \dots, x_m$.
For each $i \in \{ 1,\ldots, m\}$ there is a word $w_i=w_i(x_1, \ldots, x_m) \in [F,F]$ such that $g_i=w_i(g_1,\ldots, g_m)$.
In this case $\G_{d,k-1}$ is a homomorphic image of the perfect group $L$ with presentation \[ L:=\langle x_1, \ldots, x_m \  | \  x_1^{-1}w_1, \ldots , x_m^{-1}w_m \rangle. \] Hence we may replace $\G_{d,k-1}$ with $L$ and $g_i$ with the image of $x_i$ in $L$ proving the claim. 

Further by replacing $L$ with a free product $L*L * \cdots *L$ we may assume that $m \geq d$.
For each $j$ we pick a surjective homomorphism $f_j: L \rightarrow G_j/W_j$ and let $a_{i,j} \in G_j$ be such that $a_{i,j}W_j=f_j(g_i)$. 
For every $i\in \{1, \ldots, m\}$ and $j \in \mathcal Y$ the element $a_{i,j}w_i(a_{1,j}, \ldots , a_{m,j})^{-1}$ belongs to $W_j=A_j \times S_j$. We define $k_{i,j} \in A_j$ and $s_{i,j} \in S_j$ so that

\[ a_{i,j}w_i(a_{1,j}, \ldots , a_{m,j})^{-1}= k_{i,j} s_{i,j}.\]

Let $B_j=[A_j,G_j]$ and note that $A_j/B_j$ is central in the perfect group $G_j/B_j$. Therefore $A_j/B_j \leq \Phi(G_j/B_j)$, the Frattini subgroup of $G_j/B_j$.

Note that since $w_i \in [F,F]$ we have 
\[ w_i(g_1, \ldots, g_m) \equiv w_i(g_1c_1, \ldots, g_mc_m) \quad \textrm{mod } B_j\]
for all $g_1,\ldots ,g_m \in G_j$ and $c_1, \ldots, c_m \in A_j$.
Therefore may replace $a_{i,j}$ with $a_{i,j} k_{i,j}^{-1}$ and assume from now on that $k_{i,j} \in B_j$ for every $i,j$. Moreover $\{a_{1,j}, \ldots, a_{m,j}\}$ generates $G_j$ modulo $B_jW_j$.  Since $m \geq d$ the Gaschutz Lemma allows us to  adjust each $a_{i,j}$ by an element from $B_jW_j$ and further assume that $\langle a_{1,j}, \ldots, a_{m,j} \rangle =G_j$ for all $j \in \mathcal Y$.

Let $\mathfrak G= \prod_{j \in \N} G_j$ and let $\mathbf {a_i}= (a_{i,j})_{j \in \N} \in \mathfrak G$. Let $\Delta= \langle \mathbf{a_1}, \ldots, \mathbf{a_m}\rangle < \mathfrak G$.
We will find a finitely generated perfect group $\G < \mathfrak G$ which contains $\Delta$. Since $\Delta$ maps onto each factor $G_j$ of $\mathfrak G$ so does $\G$ and we may take $\Gamma_{d,k}=\G$.

Let $\mathfrak A= \prod_{j \in \N} A_j$, $\mathfrak S= \prod_{y \in \N} S_j$ with elements $\mathbf{k_i}:= (k_{i,j})_{j \in \N} \in \mathfrak A$ and $\mathbf{s_i}=(s_{i,j})_{j \in \N} \in \mathfrak S$.
Note that 

\begin{equation} \label{perfect}
\mathbf{k_i} \mathbf {s_i}= \mathbf{a_i}w_i( \mathbf{a_1}, \ldots, \mathbf{a_m})^{-1},  \quad i=1, \ldots, m. \end{equation}

We will find 

1. A finitely generated $\mathbb Z \Delta$-submodule $Q$ in $\mathfrak{A}$ such that $Q=[Q,\Delta]$ and $\mathbf{k_1}, \ldots, \mathbf{k_m} \in Q$, and

2. A finitely generated perfect group $T < \mathfrak S$ containing $\mathbf{s_1}, \ldots, \mathbf{s_m}$.

We take $\G$ to be the subgroup generated by $\Delta \cup  Q \cup T$ in $\mathfrak G$.
To see that $\G$ is a perfect group observe that since $T$ is perfect we have $T < [\G, \G]$. We have $Q=[Q, \Delta] = [Q, \G] <[\G, \G]$. Finally from (\ref{perfect}) and $\mathbf{k_i}\mathbf{s_i} \in QT < [\G, \G]$ we deduce that each $\mathbf{a_i} \in [\G, \G]$. Therefore $\Delta < [\G, \G]$ and hence $[\G, \G]$ contains $\Delta, Q$ and $T$. It follows that $[\G, \G]=\G$ as required. 

To find the module $Q < \mathfrak A$ we recall that each $k_{i,j} \in B_j=[A_j, G_j]$. Hence 

\[ k_{i,j}= \prod_{l=1}^m{ [q_{i,j,l}, a_{l,j}]}\] for some $q_{i,j,l} \in A_j$.
Let $\mathbf{q_{i,l}}=(q_{i,j,l})_{j \in \mathcal Y} \in \mathfrak A$ and let $Q_1$ be the $\Z \Delta$-submodule of $\mathfrak A$ generated by all $\mathbf{q_{i,l}}$ with $i,l \in \{1, \ldots, m\}$. Put $Q=[Q_1, \Delta]$ and note that 
\[\mathbf{k_i}= \prod_{l=1}^m [\mathbf {q}_{i,l}, \mathbf{a}_l] \in Q. \]

Note that the inductive hypothesis implies that the image $\bar \Delta= \Delta \mathfrak A \mathfrak S/\mathfrak A \mathfrak S$ of $\Delta$ in $\mathfrak G/\mathfrak A \mathfrak S$ is a perfect group. Since $\mathfrak{A}$ is $\Z \bar \Delta$-module it follows from Lemma \ref{G-1} (c) that $Q=[Q_1, \Delta]=Q_1(\bar \Delta-1)$ is a perfect $\Z \bar \Delta$-module and therefore $Q=[Q,\Delta]$ as required. Finally $Q$ is generated as a $\Z \Delta$-module by the elements $\mathbf{q}_{i,l}(\mathbf a_t-1)=[\mathbf q_{i,l},\mathbf a_t]$ for all $i,l,t \in \{1, \ldots, m\}$.

We now prove the existence of the subgroup $T$.
Changing notation write $\mathfrak S= \prod_{i \in \N} M_i$ as a product of nonabelian finite simple groups $M_i$, possibly with repetitions. For $l=1 \ldots, m$ let $\mathbf s_l=(s_{l,i})_{i \in \N} \in \mathfrak S$ with $s_{l,i} \in M_i$.

Let $P$ be the perfect group in the statement of Theorem \ref{seg} with generators $p_1, \ldots, p_{61}$. For each $K_i$ choose a surjective homomorphism $h_i: P \rightarrow M_i$ and let $m_{i,j}=h_i(p_j)$ for $j=1, \ldots, 61$. Let $\mathbf m_j= (m_{i,j})_{ i \in \N} \in \mathfrak S$ and let $W=\langle \mathbf m_1, \ldots, \mathbf m_{61} \rangle$. Thus $W< \mathfrak S$ is a homomorphic image of $P$ which maps onto each simple factor $M_i$.

For each $i \in \N$ the elements $\{m_{i,1}, \ldots, m_{i,61}\}$ are a generating set of $M_i$ and hence there exists some $m_{i,j_i}$ such that $|C_{M_i}(m_{i,j_i})| \leq |M_i|^{60/61}$. It follows that the conjugacy class $m_{i,j_i}^{M_i}$ of $m_{i,j_i}$ in $M_i$ has size at least $|M_i|^{1/61}$. Theorem \ref{cover} gives that there exists a constant $e \in \N$ such that the covering number of $M_i$  with respect to the set 
\[ \prod_{j=1}^{61} m_{i,j}^{M_i} \]
is at most $e$.
In particular for every $l=1, \ldots, m$ and $ i \in \N$ there exist elements $r_{l,i,j,t} \in M_i$ ($1 \leq j \leq 61 , 1 \leq t \leq e$) such that
\[ s_{l,i}= \prod_{t=1}^ e  \prod_{j=1}^{61} m_{i,j}^{r_{l,i,j, t}} \]

Let $\mathbf r_{l,j,t}: =(r_{l,i,j,t})_{i \in \N} \in \mathfrak S$ and let $W_{l,j,t}=W^{\mathbf r_{l,j,t}}$, a conjugate of $W$ in $\mathfrak S$.
Note that each $W_{l,j,t}$ is a finitely generated perfect group.

We have \[ \mathbf s_l= \prod_{t=1}^e \prod_{j=1}^{61} \mathbf m_j^{\mathbf r_{l,j,t}} \in \prod_{t=1}^e \prod_{j=1}^{61} W_{l,j,t}.\]

Therefore we define $T$ to be the subgroup generated by all $W_{l,j,t}$ for $1 \leq l \leq m, 1 \leq j \leq 61$ and $1 \leq t \leq e$.

\end{document}